\documentclass[12pt,letterpaper,oneside]{amsart}
\usepackage{euscript,verbatim,fullpage,amscd}

\newenvironment{pf}{\proof[\proofname]}{\endproof}
\newenvironment{pf*}[1]{\proof[#1]}{\endproof}

\newtheorem{theorem}{Theorem}[section]

\newtheorem{lemma}[theorem]{Lemma}
\newtheorem{conjecture}[theorem]{Conjecture}
\newtheorem{algorithm}[theorem]{Algorithm}

\theoremstyle{definition}
\newtheorem{example}[theorem]{Example}

\theoremstyle{remark}
\newtheorem{remark}{Remark}

\newcommand{\B}{\operatorname{B}}
\newcommand{\Zar}{\operatorname{Zar}}
\newcommand{\characteristic}{\operatorname{char}}

\newcommand{\orb}{\operatorname{orb}}
\newcommand{\im}{\operatorname{im}}
\newcommand{\rank}{\operatorname{rank}}

\newcommand{\Hom}{\operatorname{Hom}}
\newcommand{\Pic}{\operatorname{Pic}}
\newcommand{\Cl}{\operatorname{Cl}}
\newcommand{\et}{\operatorname{\acute{e}t}}
\newcommand{\SF}{\operatorname{SF}}
\newcommand{\emb}{\operatorname{emb}}
\newcommand{\Cech}{\operatorname{\check Cech}}
\newcommand{\scr}{\EuScript}

\begin{document}
\def \TNemb {T_N\emb}
\def \dis {\displaystyle}
\def \scrSF {\scr{S}\scr{F}}
\def \scrW {\scr{W}}
\def \scrO {\scr{O}}
\def \scrP {\scr{P}}
\def \scrU {\scr{U}}
\def \scrM {\scr{M}}

\title[Invariants of a Fan]{Topological invariants of a fan associated to a
toric variety}
\author[T. J. Ford]{Timothy J. Ford}

\address{Department of Mathematics, Florida Atlantic University,
Boca Raton, Florida 33431}
\email{Ford@@acc.fau.edu}

\begin{abstract}
Associated to a toric variety $X$ of dimension $r$ over a field $k$ is a
fan $\Delta$ on
$\Bbb R^r$. The fan $\Delta$ is a finite set of cones which are in
one-to-one correspondence with the orbits of the torus action on $X$. The
fan $\Delta$ inherits the Zariski topology from $X$. In this article some
cohomological invariants of $X$ are studied in terms of whether or not they
depend only on $\Delta$ and not $k$. Secondly some
numerical invariants of $X$ are studied in terms of whether or not they are
topological invariants of the fan $\Delta$. That is, whether or not they
depend only on the finite topological space defined on $\Delta$. The
invariants with which we are mostly concerned are
the class group of Weil divisors, the Picard group, the Brauer group and
the dimensions of the
torsion free part of the \'etale cohomology groups with coefficients in the
sheaf of units.  The notion of an open neighborhood of a fan is introduced
and examples are given for which the above invariants are sufficiently fine
to give nontrivial stratifications of an open neighborhood of a fan all of
whose maximal cones are nonsimplicial.
\end{abstract}
\maketitle

\section{Introduction} \label{sec1}
Let $k$ be a field.
Let $N = \Bbb Z^r$ and denote by $T_N$ the $k$-torus on $N$. Let $\Delta$
be a finite fan on $N \otimes \Bbb R$ and $X = \TNemb(\Delta,k)$ the toric
variety over $k$ associated to $\Delta$ \cite{D:Gtv}, \cite{F:ITV},
\cite{O:CBA}.
This defines a functor $\TNemb$ on the product category
\begin{equation} \label{eq33}
\begin{array}{ccc}
(\text{finite fans on }N \otimes \Bbb R) \times (\text{fields}) &
\stackrel{\TNemb}{\longrightarrow} & (\text{toric varieties}) \\
(\Delta,k)  & \mapsto & \TNemb(\Delta,k)
\end{array}
\text{.}
\end{equation}
We define the topology on $\Delta$
as follows (cf. \cite[pp. 137--138]{DFM:CBg}).
The orbit space $\tilde X$ of $X$ under the action of the torus $T_N$ is in
one-to-one correspondence with the finite set of cones that belong to
$\Delta$. There is a topology on $\tilde X$ inherited from $X$ by the
continuous function $X \to \tilde{X}$. Identifying a cone $\sigma \in
\Delta$ with the orbit $\orb{\sigma}$ in $\tilde X$, we see that the
topology on $\tilde X$ corresponds to the topology on $\Delta$ under which
the open sets are the subfans of $\Delta$. The fan $\Delta$ is now a
two-faced beast. On the one hand $\Delta$ is an object in the category of
fans on $N \otimes \Bbb R$. At the same time $\Delta$ is an object of the
category of finite topological spaces. To distinguish between these roles
played by $\Delta$, we denote by $\Delta_{fan}$ the object in the category
of fans on $N \otimes \Bbb R$ and by $\Delta_{top}$ the object in the
category of finite topological spaces. This defines a functor $\frak T$
(which
factors via $\TNemb$ through the category of toric varieties)
\begin{equation} \label{eq26}
\begin{array}{ccc}
(\text{finite fans on }N \otimes \Bbb R)  \times (\text{fields}) &
\xrightarrow{\frak T} &
(\text{finite top. spaces}) \\
(\Delta_{fan},k)  & {\mapsto} & \Delta_{top}
\end{array} \text{.}
\end{equation}

In Section~\ref{sec2} we consider some invariants of $\Delta_{fan}$ that
are constant for all $k$.
Suppose $\gamma(\Delta_{fan},k)$
is an invariant
that is defined for any pair  $(\Delta_{fan},k)$ (in this article $\gamma$
is usually an abelian group). We call $\gamma$ a {\em fan invariant} in
case  $\gamma(\Delta_{fan},k)$ depends only on $\Delta_{fan}$ and not on
$k$ --- that is, given a fan  $\Delta_{fan}$,  $\gamma(\Delta_{fan},k_1)
\cong  \gamma(\Delta_{fan},k_2)$ for every pair of fields $k_1$, $k_2$.
We show that the Brauer group
$\B(~)$ is not a fan invariant for nonsingular fans. This is an observation
based on a theorem of Hoobler and \cite[Theorem~1.1]{DF:Bgt}.
In \cite[Theorem~1.1]{DF:Bgt} a complete computation of the Brauer group of
a nonsingular toric variety $X=\TNemb(\Delta)$ over an algebraically closed
field $k$ was given in terms of the
so-called invariant factors of the fan $\Delta$. In Theorem~\ref{th5} we
give the Brauer group of $\TNemb(\Delta,k)$ for any field $k$ in terms of
the Brauer group and Galois group of $k$.
The main result of
Section~\ref{sec2} is Theorem~\ref{th6} in which it is stated
that the class group, $\Cl(~)$, the Picard group, $\Pic(~)$, and the
relative cohomology group, $H^2(K(~)/X(~)_{\et},\Bbb G_m)$ (where $K(X)$
is the function field of $X$), are fan
invariants. The proof of Theorem~\ref{th6} follows from that of
\cite[Theorem~1]{DFM:CBg} and is omitted.

In Section~\ref{sec3} we consider some invariants of $\Delta_{fan}$ that are
constant on fibers of the map $\frak T$ in \eqref{eq26}, hence
depend only
on $\Delta_{top}$. That is, suppose we have an invariant
$\beta(\Delta_{fan})$ (usually a numerical invariant) associated to any fan
$\Delta_{fan}$. If two fans $\Delta_1$, $\Delta_2$ have the same
$\beta$-invariant whenever $(\Delta_1)_{top} \cong (\Delta_2)_{top}$, then
we say $\beta$ is a {\em topological invariant of $\Delta_{fan}$}.
We consider several invariants, all being cohomologically defined. The
first sequence is defined by
\[
\rho_0 = \dim_{\Bbb Q} \left[ H^0(X_{\et},\Bbb G_m) / k^* \otimes \Bbb Q
\right] \text{,}
\]
and for $i \ge 1$,
\[
\rho_i = \dim_{\Bbb Q} \left[ H^i(X_{\et},\Bbb G_m) \otimes \Bbb Q \right]
\text{.}
\]
Also set
\[
\rho_1' = \dim_{\Bbb Q} \left[ \Cl(X) \otimes \Bbb Q  \right] \text{.}
\]
For $0 \le i \le 2$ these numbers are finite and are fan invariants.
The first main result of Section~\ref{sec3} lists some facts about $\rho_0$
and $\rho_1'$.

\medskip\noindent
{\bf Theorem~\ref{th1}.}
{\em
Let $N = \Bbb Z^r$, $\Delta$ a fan on $N \otimes \Bbb R$,
$X = \TNemb(\Delta)$, and
$s =
\dim_{\Bbb R} \Bbb R|\Delta_{fan}|$ (that is, $s$ is the dimension of the $\Bbb
R$-vector space spanned by the vectors in the support $|\Delta_{fan}|$). Then
\begin{itemize}
\item[(a)]
$\rho_0 = r-s$, hence is a fan invariant, but not a topological invariant.
\item[(b)]
Suppose $\Delta_{fan}$ contains a cone $\sigma$ such that $\dim{\sigma} =
r$. This is true for example if $\Delta_{fan}$ is a complete fan on $N
\otimes \Bbb R$. Then $\rho_0 = 0$ and $\rho_0$ is a topological invariant
of $\Delta_{fan}$.
\item[(c)]
$\rho_1'=\dim_{\Bbb Q}(\Cl(X) \otimes \Bbb Q) = \#(\Delta(1)) -s$.
If $\dim{\Delta_{top}} =r$, then  $\rho_1'$ is a topological invariant of
$\Delta_{fan}$.
\item[(d)]
The number $\rho_0-\rho_1'$ is a topological
invariant of $\Delta_{fan}$.
\end{itemize}
}

The second main result of Section~\ref{sec3} gives some results on
$\rho_0$, $\rho_1$ and $\rho_2$ for simplicial fans.

\medskip\noindent
{\bf Theorem~\ref{th2}.~}{\em
Let $N = \Bbb Z^r$.
Let $\Delta$ be a simplicial fan on $N \otimes \Bbb R$
and $s = \dim_{\Bbb R}\Bbb R |\Delta_{fan}|$.
Then
\begin{itemize}
\item[(a)]
 $\rho_1 =
\#(\Delta(1))-s$.
\item[(b)]
  $\rho_2 = 0$ hence is a topological invariant of
$\Delta_{fan}$.
\item[(c)]
If $\dim{\Delta_{top}}=r$,
then $\rho_1$ is a topological invariant of $\Delta_{fan}$.
\item[(d)]
$\rho_0-\rho_1+\rho_2$ is a topological
invariant of  $\Delta_{fan}$.
\end{itemize}}

The third main result of Section~\ref{sec3} gives some results on
$\rho_0$, $\rho_1$ and $\rho_2$ for 3-dimensional fans.

\medskip\noindent
{\bf Theorem~\ref{th3}.~}{\em
Let $\Delta$ be a fan on $N \otimes \Bbb R$.
Let $\sigma_0, \dots, \sigma_w$ be the maximal cones in $\Delta$.
Assume
$\sigma_i \cap
\sigma_j$ is simplicial for each $i \not = j$. These assumptions are
satisfied for example if $\dim{\Delta_{top}} \le 3$.
Then
\begin{itemize}
\item[(a)]
\[
 \rho_1 +s + \sum_{i=0}^w(\#(\Delta(\sigma_i)(1))-s_i) =
 \rho_2 +  \#( \Delta(1)) \text{,}
\]
where  we set $s_i = \dim{\sigma_i}$ for each $i = 0, \dots, w$ and $s =
\dim_{\Bbb R} \Bbb R |\Delta_{fan}|$.
\item[(b)]
$\rho_0-\rho_1+\rho_2$ is a topological invariant of $\Delta_{fan}$.
\end{itemize}}

In Section~\ref{sec4} we introduce the notion of an {\em open neighborhood
$B$ of a fan $\Delta$}. This is a subset of the fiber $\frak
T^{-1}(\Delta_{top})$ that is parametrized by a dense subset of a real
manifold.
Let $\scrSF$ denote the sheaf of $\Delta$-linear support functions on the
topological space $\Delta_{top}$. It was shown in
\cite{DFM:CBg} that the numbers $\rho_i$, $1\le i\le 2$,  can be
determined by the cohomology of the sheaf $\scrSF$ on the finite
topological space $\Delta_{top}$.  Therefore we define another sequence of
invariants by
\[
\kappa_i = \dim_{\Bbb Q}\left[ H^i(\Delta_{top},\scrSF) \otimes \Bbb Q
\right]
\]
for $i \ge 0$.
We consider the stratification of $B$ by the numerical invariant
$\kappa_0$.
Several examples are given for which the stratification
of $B$ is nontrivial.
We conjecture that $\kappa_0 = 3$ on a nonempty open subset of $B$ if
$\Delta$ is a complete fan on $\Bbb R^3$ such that every maximal cone of
$\Delta$  is nonsimplicial. Algorithm~\ref{alg1} is presented which
computes an upper bound for $\kappa_0$. For complete 3-dimensional fans,
this algorithm can be used to compute an upper bound for $\rho_1$ and
$\rho_2$.

For the benefit of the reader the following notation will be fixed
throughout the rest of the paper.

\begin{table}[htp]
\label{tab1}
\end{table}
\begin{center}
\begin{tabular}{|lp{14pc}|ll|}\hline
$k$ & a field          & $r$ & a positive integer \\
$N$ & $=\Bbb Z^r$                         & $M $ & $=\Hom_{\Bbb Z}(N,\Bbb Z)$
\\
$\Delta$ & \raggedright a finite rational fan on $N\otimes \Bbb R$
                                     & $X$ & $ =\TNemb(\Delta,k)$ toric variety
\\
$\Delta_{fan}$ & object in the category of fans
                                     & $\Delta_{top}$ & finite topological
space \\
$\left|\Delta_{fan}\right|$ & support of the fan $\Delta$
                  & $\frak T$ & functor that maps $(\Delta_{fan},k)$ to
$\Delta_{top}$ \\
$\Cl(X)$ &\raggedright class group of Weil divisor classes & $\Pic{X}$ & Picard
group of
invertible modules \\
$\B(X)$ & \raggedright Brauer group of Azumaya algebra classes
 & $\Bbb G_m$ & \'etale sheaf of units \\
$\rho_0 $ & $ = \dim_{\Bbb Q}\left[ H^0(X_{\et},\Bbb G_m)\otimes\Bbb Q\right]$
&
$\rho_i $ &  $ = \dim_{\Bbb Q}\left[ H^i(X_{\et},\Bbb G_m)\otimes\Bbb
Q\right]$ (for $i>0$) \\
$\rho'_1$ & $ = \dim_{\Bbb Q}\left[ \Cl(X)\otimes\Bbb Q\right]$ &
$\kappa_i$ &  $=\dim_{\Bbb Q}\left[
H^i(\Delta_{top},\scrSF)\otimes\Bbb Q\right]$
(for $i\ge 0$) \\
$s $ & $=\dim_{\Bbb R}\Bbb R \left|\Delta_{fan}\right|$ & $\Delta(i)$ &
$=\{\sigma\in\Delta | \dim{\sigma}=i\}$ \\
$K$ & $=K(X)$ the function field of $X$ & $\scrSF$ & sheaf of $\Delta$-linear
support functions \\
$\scrW$ & sheaf of Weil divisors & $\scrP$ & quotient sheaf $\scrW / \scrSF$ \\
\hline
\end{tabular}\end{center}

\section{Fan Invariants} \label{sec2}
In Theorem~\ref{th5} we determine the Brauer group of a nonsingular toric
variety over $k$. This invariant depends on $k$. We then show in
Theorem~\ref{th6} that $\Cl(X)$, $\Pic{X}$ and the relative cohomology
group $H^2(K/X_{\et},\Bbb G_m)$ depend only on $\Delta_{fan}$, not on $k$.
In order to determine the Brauer group of a nonsingular toric variety over
$k$, we use the following theorem of Hoobler.

\begin{theorem}
\label{th4}
Let $R = A[x_1,x_1^{-1}, \dots, x_r,x_r^{-1}]$,
where $A$ is a connected,
normal integral domain. Suppose $\nu$ is an integer relatively prime to  the
residue characteristics of $A$. Then
\begin{equation}
\label{eq8}
H^1(R,\Bbb Z/\nu) = H^1(A,\Bbb Z/\nu) \oplus \left( \bigoplus^r \mu_\nu^{-1}
\right) \text{ ,}
\end{equation}
and
\begin{equation}
\label{eq27}
_\nu\B(R) = {_\nu\B(A)} \oplus \left( \bigoplus^r H^1(A, \Bbb Z/\nu)
\right) \oplus \left( \bigoplus^{r(r-1)/2} \mu_\nu^{-1} \right) \text{ .}
\end{equation}
\end{theorem}
\begin{pf} See \cite[Cor. 2.6]{H:Fgr}. \end{pf}

Therefore the $\nu$-torsion  of the Brauer  group  of $R=A[x_1,x_1^{-1},
\dots, x_r,x_r^{-1}]$    is
generated by the Azumaya $A$-algebras and the classes of cyclic crossed
product algebras of 2 types. For each cyclic Galois extension $C/A$ of
degree $\nu$ with group $\langle \sigma \rangle$ and for each $1 \le i \le
r$, there is the cyclic crossed
product $(C/R,\langle \sigma \rangle, x_i)$ which is an Azumaya algebra
over $R$. If there exists a primitive $\nu$-th root of unity $\zeta$ over
$A$, then the symbol algebras $(x_i,x_j)_\nu$ are Azumaya algebras over $R$.

\begin{example}
\label{ex13}
Let  $R= \Bbb R [x_1,x_1^{-1},\dots, x_r,x_r^{-1}]$. Then by
Theorem~\ref{th4}
$\B(R)$ is an elementary 2-group and
\begin{equation}
\label{eq28}
\begin{split}
\B(R) & \cong {\B(\Bbb R)} \oplus \left( \bigoplus^r H^1(\Bbb R, \Bbb Z/2)
\right) \oplus \left( \bigoplus^{r(r-1)/2} \mu_2^{-1} \right) \\
 & \cong \left( \Bbb Z/2 \right)^{1+r+
r(r-1)/2}
\end{split}
\end{equation}
\end{example}

Define a sheaf ${\scrSF}$  on $\Delta_{top} $
by assigning to each open set $\Delta ' \subseteq \Delta_{top}$
the abelian group $\SF(\Delta ')$ of support functions on $\Delta '$.
Let $M = \Hom(N,{\Bbb Z})$ be the dual of $N$.
There is a natural map $M \rightarrow \SF(\Delta ')$
which is locally surjective.
If ${\scrM}$  denotes the constant sheaf of $M$ on $\Delta_{top}$,
then there is an exact sequence of sheaves on $\Delta_{top}$:
\begin{equation}
\label{eq32}
0 \to  {\scrU}  \to  {\scrM}  \to {\scrSF}  \to  0
\end{equation}
\noindent
where ${\scrU}$ is defined by the sequence \eqref{eq32}.
On any open $\Delta ' \subseteq  \Delta_{top}$,
${\scrU} (\Delta ') = |\Delta '|^\perp  \cap  M =
\{m \in  M | \langle m,y\rangle  = 0$ for all $y \in  |\Delta '|\}$.
Because ${\scrM}$  is flasque,
$H^p(\Delta_{top} ,{\scrM} ) = 0$ for all $p \geq  1$, so
$H^p(\Delta_{top} ,{\scrSF} ) \cong  H^{p+1}(\Delta_{top} ,{\scrU} )$
for all $p \geq  1.$

Let $k$ be a field and $X = \TNemb(\Delta,k)$ a nonsingular toric variety
over $k$.
Let
$N' = \langle \bigcup_{\sigma \in \Delta} \sigma \cap N \rangle$, let $\nu
\ge 2$ be relatively prime to $\characteristic{k}$,
and let $M_\nu = \{ m \in M | \langle m,n' \rangle \equiv 0 \pmod \nu
\text{ for all } n' \in N' \}$.
The basis theorem for finitely generated abelian groups gives a basis $n_1,
\dots, n_r$ of $N$ such that
$N' = \Bbb Z a_1 n_1 \oplus \Bbb Z a_2 n_2 \oplus \dots \oplus \Bbb Z a_r n_r$
where the $a_i$ are nonnegative integers and
$a_i | a_{i+1}$ for $1 \le i \le r-1$.
As in \cite{DF:Bgt} call  $a_1, \dots, a_r$  the set of invariant factors
of $X$.

\begin{theorem} \label{th5} In the above terminology,
if $(\nu, a_i)$ is the greatest common
divisor of $\nu$ and $a_i$, then
\begin{equation}
\label{eq29}
\begin{split}
H^1(X,\Bbb Z/\nu) &\cong
H^1(k,\Bbb Z/\nu) \oplus
\left(
M_\nu/\nu M \otimes \mu_\nu  ^{-1}
\right) \\
&\cong
H^1(k,\Bbb Z/\nu) \oplus \left(
\bigoplus_{i=1}^r \Bbb Z/(\nu,a_i) \otimes \mu_\nu ^{-1}
\right) \text{ .}
\end{split}
\end{equation}
\begin{equation}
\label{eq30}
\begin{split}
{_\nu \B(X)} & = {_\nu \B'(X)} \cong \\ &
{_\nu \B(k)} \oplus
\left( \bigoplus_{i=1}^r H^1(k,\Bbb Z/\nu) \otimes \Bbb Z/(\nu,a_i)
\right) \oplus
\left(
\bigoplus_{i=1}^r \Hom(\Bbb Z/a_i \otimes \mu_\nu, \Bbb Q/ \Bbb Z )^{r-i}
\right)
\end{split}
\end{equation}
\end{theorem}
\begin{pf} Follows from the proof of \cite[Theorem~1.1]{DF:Bgt} and
Theorem~\ref{th4}. \end{pf}

Therefore the $\nu$-torsion  of the Brauer  group  of the nonsingular toric
variety $X$ is
generated by the classes of algebras from $k$ and
cyclic crossed
product algebras of 2 types. For each cyclic Galois extension $C/k$ of
degree $\nu$ with group $\langle \sigma \rangle$ and for each $1 \le i \le
r$, there is the cyclic crossed
product $(C/k,\langle \sigma \rangle, x_i)$ which is an Azumaya algebra
over the torus $T_N$. This algebra is unramified on $X$ if and only if the
function $x_i$ corresponds to an element of $M_\nu$.
If there exists a  $\nu$-th root of unity $\zeta$ over
$k$, then the symbol algebras $(x_i,x_j)_\nu$ are Azumaya algebras over
$T_N$. Those symbols which  are unramified on $X$  correspond to the
last summand of \eqref{eq30}.

\begin{example}
\label{ex14}
Let  $k= \Bbb R$ and $X = \TNemb(\Delta)$ a nonsingular toric variety over
$\Bbb R$.
Then by
Theorem~\ref{th5}
$\B(X)$ is an elementary 2-group. If $t = \left| \{ a_i | (2,a_i) \not = 1
\} \right|$, then
\begin{equation}
\label{eq31}
\begin{split}
 \B(X) & \cong
 \B(\Bbb R) \oplus
\left( \bigoplus_{i=1}^r  \Bbb Z/(2,a_i)
\right) \oplus
\left(
\bigoplus_{i=1}^r \Hom(\Bbb Z/a_i \otimes \mu_2, \Bbb Q/ \Bbb Z )^{r-i}
\right) \\
 & \cong \Bbb Z/2 \oplus (\Bbb Z/2)^t \oplus (\Bbb Z/2)^{t(t-1)/2}
\end{split}
\end{equation}
\end{example}

\begin{theorem} \label{th6}
Let $k$ be a field and $X = \TNemb(\Delta)$ a toric variety
over $k$ with function field $K$.
Then
\begin{enumerate}
\item
$H^p(\Delta_{top},\scrU) \cong H^p(X_{\Zar},\scrO^*)$ for all $p \ge 1$ hence
$H^p(X_{\Zar},\scrO^*)$ depends only on $\Delta_{fan}$, not $k$. In particular
$\Cl(X)$ and $\Pic{X}$ depend only on $\Delta_{fan}$.
\item
$H^1(\Delta_{top},\scrSF) \cong H^2(X_{\Zar},\scrO^*) \cong H^2(K/X_{\et},\Bbb
G_m)$ hence $H^2(K/X_{\et},\Bbb G_m)$ depends only on $\Delta_{fan}$, not $k$.
\item
If  $\tilde \Delta $ is a nonsingular subdivision of $\Delta $
and $\tilde X = \TNemb(\tilde \Delta )$, then the sequence
\[
0 \to  H^2(K/X_{\et},\Bbb G_m) \to
H^2(X_{\et},\Bbb G_m) \to
H^2(\tilde X_{\et},\Bbb G_m) \to  0
\]
(with natural maps) is split-exact.
\end{enumerate}
\end{theorem}
\begin{pf}
The theorem follows from
\cite{DFM:CBg},
noting that the proof of \cite[Theorem~1]{DFM:CBg} did not assume
that $k$ is algebraically closed until the proof of Lemma~7 where it was
not necessary anyway.
\end{pf}

\section{Topological Invariants} \label{sec3}
The first invariants to be considered as candidates for topological
invariants are the following. Let $\Delta$ and $X$ be as in the
Introduction. For each $i \ge 0$ we define a positive
integer $\rho_i$.
Set
\[
\rho_0 = \dim_{\Bbb Q} \left[ H^0(X_{\et},\Bbb G_m) / k^* \otimes \Bbb Q
\right] \text{,}
\]
\[
\rho_1 = \dim_{\Bbb Q} \left[ H^1(X_{\et},\Bbb G_m) \otimes \Bbb Q \right]
= \dim_{\Bbb Q} \left( \Pic{X} \otimes \Bbb Q \right)
\]
and for $i \ge 2$,
\[
\rho_i = \dim_{\Bbb Q} \left[ H^i(X_{\et},\Bbb G_m) \otimes \Bbb Q \right]
\text{.}
\]
The number $\rho_1$ is the traditional Picard number $\rho$ associated to
$X$. Also set
\[
\rho_1' = \dim_{\Bbb Q} \left[ \Cl(X) \otimes \Bbb Q  \right] \text{.}
\]

It follows from Theorem~\ref{th1} below that $\rho_0$ is a fan invariant
and from Theorem~\ref{th6} above that $\rho_1$, $\rho_1'$, and $\rho_2$ are
fan invariants.
Since $\Delta$ is finite, $\rho_0$, $\rho_1$, $\rho_1'$ and $\rho_2$ are
finite. For $\rho_0$, $\rho_1$ and $\rho_1'$ see \cite{O:CBA} or
\cite{F:ITV}. For $\rho_2$ this follows from \cite{DFM:CBg}.

Examples where the number $\rho_2$ is computed seem to be somewhat scarce.
Grothendieck \cite[II]{G:GB} and Childs \cite{C:Bgn} each give an example
of a local ring $\scrO_x$ on a normal surface where
$H^2((\scrO_x)_{\et},\Bbb G_m)$
is torsion free, but in each case $H^2((\scrO_x)_{\et},\Bbb G_m)$ is not
finitely generated.

\begin{remark}\label{re6}
The dimension of the topological space $\Delta_{top}$ is defined to be the
length of a maximal chain of irreducible closed subsets. One can check that
this is equal to $\max{ \{ \dim{\sigma} |}$ ${ \sigma \in \Delta  \} }$.
Therefore  $\dim{\Delta_{top}}$ is a topological invariant of
$\Delta_{fan}$.
\end{remark}

\begin{remark}\label{re7}
Define another sequence of invariants by
\[
\kappa_i = \dim_{\Bbb Q}\left[ H^i(\Delta_{top},\scrSF) \otimes \Bbb Q
\right]
\]
for $i \ge 0$. It follows from Theorem~\ref{th6}~(2) that $\kappa_1 =
\rho_2$.
Let $\sigma_0, \dots, \sigma_m$ be the maximal cones in
$\Delta$. From \cite[Lemma 8]{DFM:CBg} $\kappa_i$ can be computed from the
$\Cech$
complex
\begin{equation}
\label{eq12}
0 \to \underset{i}{\oplus}
\scrSF(\sigma_i) \stackrel{\delta^0}{\rightarrow} \underset{i<j}{\oplus}
\scrSF(\sigma_{ij}) \stackrel{\delta^1}{\rightarrow} \underset{i<j<k}{\oplus}
\scrSF(\sigma_{ijk}) \to \dots
\end{equation}
For any cone $\tau \in \Delta$, $\dim_{\Bbb Q}(\scrSF(\Delta(\tau)) \otimes
\Bbb Q) = \dim{\tau}$. Therefore, if $C^i$ denotes the $i$-th group of $\Cech$
cochains in \eqref{eq12} and $c_i = \dim_{\Bbb Q}(C^i \otimes \Bbb Q)$,
then the integer
\begin{equation}
\label{eq24}
c_0 - c_1 + c_2 - \dots
\end{equation}
is a topological invariant of $\Delta_{fan}$. Note that there exists an
integer $M$ such that $C^j = 0$ for all $j > M$. If $\dim(\Delta_{top}) =
t$, then $\kappa_j = 0$ for all $j>t$. So the left hand side of
\begin{equation}
\label{eq25}
\kappa_0 - \kappa_1 + \dots (-1)^t \kappa_t =
c_0 - c_1 + c_2 - \dots (-1)^M c_M
\end{equation}
is a topological invariant of $\Delta_{fan}$.
\end{remark}

\begin{remark}\label{re5}
Let  $\Delta$
be a finite fan on $N \otimes \Bbb R$ where $N = \Bbb Z^r$.
Setting $s = \dim_{\Bbb R} \Bbb R|\Delta_{fan}| $, we see
that $s$ is not a topological invariant of $\Delta_{fan}$. Since $s \ge
\dim{\Delta_{top}}$, if $\dim{\Delta_{top}} = r$, then $s=r$ so if $\Delta$
contains a cone $\sigma$ such that $\dim{\sigma}=r$, then $s=r$ and $s $ is
a topological invariant of $\Delta_{fan}$. This condition is satisfied, for
instance, if $\Delta$ is a complete fan on $N \otimes \Bbb R$.
\end{remark}

\begin{remark}\label{re2}
Let $\sigma \in \Delta$ and let $\Delta(\sigma)$ denote the subfan of
$\Delta$ consisting of the cone $\sigma$ and all of its faces. Then
\[
\dim{\sigma} = \dim{\Delta(\sigma)_{top}} \text{,}
\]
so the dimensions of the cones in $\Delta$ depend only on $\Delta_{top}$.
In particular, the number of 1-dimensional cones in $\Delta_{fan}$ is a
topological invariant.
\end{remark}
\begin{remark}\label{re3}
The fan $\Delta_{fan}$ is complete by definition if $|\Delta_{fan}| = \Bbb
R^r$. This is true if and only if
\begin{itemize}
\item[(i)] $\Delta(r) \not = \emptyset$ and
\item[(ii)]
for each cone $\sigma \in \Delta(r)$ and every $r-1$-dimensional face
$\tau$ of $\sigma$ there is a cone $\sigma_1 \in \Delta(r)$ such that $\tau
= \sigma \cap \sigma_1$.
\end{itemize}
But these two conditions depend only on
$\Delta_{top}$. That is, completeness can be thought of as a topological
property of $\Delta_{fan}$.
\end{remark}

{F}rom the next theorem, which combines some results on $\rho_0$ and
$\rho_1'$, we see that $\rho_0$ depends only on the dimension of
the subspace spanned by $|\Delta_{fan}|$.

\begin{theorem}
\label{th1}
Let $N = \Bbb Z^r$, $\Delta$ a fan on $N \otimes \Bbb R$,
$X = \TNemb(\Delta)$, and
$s =
\dim_{\Bbb R} \Bbb R|\Delta_{fan}|$ (that is, $s$ is the dimension of the $\Bbb
R$-vector space spanned by the vectors in the support $|\Delta_{fan}|$). Then
\begin{itemize}
\item[(a)]
$\rho_0 = r-s$, hence is a fan invariant, but not a topological invariant.
\item[(b)]
Suppose $\Delta_{fan}$ contains a cone $\sigma$ such that $\dim{\sigma} =
r$. This is true for example if $\Delta_{fan}$ is a complete fan on $N
\otimes \Bbb R$. Then $\rho_0 = 0$ and $\rho_0$ is a topological invariant
of $\Delta_{fan}$.
\item[(c)]
$\rho_1'=\dim_{\Bbb Q}(\Cl(X) \otimes \Bbb Q) = \#(\Delta(1)) -s$.
If $\dim{\Delta_{top}} =r$, then  $\rho_1'$ is a topological invariant of
$\Delta_{fan}$.
\item[(d)]
The number $\rho_0-\rho_1'$ is a topological
invariant of $\Delta_{fan}$.
\end{itemize}
\end{theorem}
\begin{pf} (a)
Let $N_1 = N \cap \Bbb R |\Delta_{fan}|$, $M_1 = \Hom_{\Bbb Z}(N_1,\Bbb
Z)$, $M_2 = N_1^{\perp}$, $N_2 = \Hom_{\Bbb Z}(M_2,\Bbb Z)$. Then $M = M_1
\oplus M_2$ and $N = N_1 \oplus N_2$. Viewing $\Delta$ as a fan on the
$s$-dimensional vector space $N_1 \otimes \Bbb R$, $X = T_{N_1}\emb(\Delta)
\times T_{N_2}$ and $H^0(X,\Bbb G_m) / k^* \cong H^0(T_{N_2},\Bbb G_m) /
k^* \cong \Bbb Z^{r-s}$. So $\rho_0 = r-s$.

(b)
In this case $s = \dim_{\Bbb R} \Bbb R|\Delta_{fan}|
= \dim_{\Bbb R} (\Bbb R \sigma)  = \dim \sigma = r$.

(c)
Let $N_0 = N \cap \Bbb R |\Delta_{fan}|$ be the set of lattice points in the
subspace $\Bbb R |\Delta_{fan}|$
and $M_0 = \Hom_{\Bbb Z}(N_0,\Bbb Z)$. Then
$\dim_{\Bbb Q}(M_0 \otimes \Bbb Q) =
\dim_{\Bbb Q}(N_0 \otimes \Bbb Q) = \dim_{\Bbb R} \Bbb R |\Delta_{fan}| =
 s$. From \cite[Corollary~2.5]{O:CBA} there is a presentation
of $\Cl(X)$
\begin{equation}
\label{eq0}
0 \to M_0 \to \bigoplus_{\rho \in \Delta(1)} \Bbb Z \rho \to
\Cl(X) \to 0 \text{.}
\end{equation}
So
$\dim_{\Bbb Q}(\Cl(X) \otimes \Bbb Q)
 = \#(\Delta(1)) -s$. In
particular, if $\dim{\Delta_{top}} =r$, then $r=s$ so  $\dim_{\Bbb Q}(\Cl(X)
\otimes \Bbb Q)$ is a topological invariant.

(d)
{F}rom (a) and (c),
\[
\rho_0-\rho_1' = (r - s) - (\#(\Delta(1)) -s) = r  - \#(\Delta(1))
\]
which depends only on $\Delta_{top}$.
\end{pf}

\begin{remark}
\label{re1}
Let $N = \Bbb Z^r$, $\Delta$ a fan on $N \otimes \Bbb R$, $s = \dim_{\Bbb
R} \Bbb R |\Delta_{fan}|$, $t = \underset{\sigma \in \Delta}{\max} \{
\dim{\sigma} \}$.
It follows from \cite[Theorem~2.3]{DF:Bgt} that if
$t \le 2$, then $\rho_1 = \#(\Delta(1)) - s$ and $\rho_2 = 0$.
In this case $\Delta$ is a simplicial fan, so this is a
special case of the following theorem.
\end{remark}

\begin{theorem}
\label{th2}
Let $N = \Bbb Z^r$.
Let $\Delta$ be a simplicial fan on $N \otimes \Bbb R$
and $s = \dim_{\Bbb R}\Bbb R |\Delta_{fan}|$.
Then
\begin{itemize}
\item[(a)]
 $\rho_1 =
\#(\Delta(1))-s$.
\item[(b)]
  $\rho_2 = 0$ hence is a topological invariant of
$\Delta_{fan}$.
\item[(c)]
If $\dim{\Delta_{top}}=r$,
then $\rho_1$ is a topological invariant of $\Delta_{fan}$.
\item[(d)]
$\rho_0-\rho_1+\rho_2$ is a topological
invariant of  $\Delta_{fan}$.
\end{itemize}
\end{theorem}
\begin{pf}
Suppose $\sigma \in \Delta$ is a simplicial cone and
$\dim{\sigma} = d$. For any support function $h \in \scrSF(\Delta)$,
$h|_\sigma$ is linear and completely determined by its values on  a
spanning set $\{\eta_1, \dots, \eta_d \} \subseteq N$ for $\sigma$. Since
$\dim{\sigma}=d$, $\sigma$ is spanned by $d$ lattice points. So
$\scrSF(\Delta(\sigma)) \otimes \Bbb Q \cong \Bbb Q^d$.

If $\tau_0, \dots, \tau_n$ are the cones in $\Delta(1)$, and $\Gamma =
\{0,\tau_0, \dots, \tau_n \}$, then $\Gamma_{top}$ is an open subset of the
topological space $\Delta_{top}$. Define the sheaf $\scrW$ on
$\Delta_{top}$ to be the direct image $i_*(\scrSF |_{\Gamma_{top}})$.
Since $\Gamma_{fan}$ is a nonsingular fan, $\scrSF |_{\Gamma_{top}}$ is the
sheaf defined by $\Xi \mapsto \Bbb Z^{\#(\Xi(1))}$ for each open $\Xi
\subseteq \Delta_{top}$.
It follows that $\scrW(\Xi) = \Bbb Z^{\# (\Xi(1))}$, hence $\scrW$ is a
flasque sheaf. So there is an embedding $\scrSF \to \scrW$ of sheaves on
$\Delta_{top}$ and we define $\scrP$ by the exact sequence of sheaves
\cite[(13), p. 149]{DFM:CBg}
\begin{equation}
\label{eq1}
0 \to \scrSF \to \scrW \to \scrP \to 0  \text{.}
\end{equation}

Since $\Delta$ is simplicial,
$ \scrSF(\Delta(\sigma))$ and $\scrW(\Delta(\sigma))$
are  free of the same rank $\dim{\sigma}$. Therefore, $\scrP$ is
locally torsion, hence torsion.
Because $\scrW$ is flasque, $H^1(\Delta_{top},\scrW) = 0$ and the long exact
sequence associated to \eqref{eq1} becomes
\begin{equation}
\label{eq2}
0 \to H^0(\Delta_{top},\scrSF) \to H^0(\Delta_{top},\scrW) \to
H^0(\Delta_{top},\scrP)
\to H^1(\Delta_{top}, \scrSF) \to 0  \text{.}
\end{equation}
Because $\scrP$ is torsion, $H^0(\Delta_{top},\scrP)$ is torsion. So
$H^1(\Delta_{top}, \scrSF) \otimes \Bbb Q = 0$. By
\cite[Theorem~1]{DFM:CBg} $\rho_2 = \dim(H^1(\Delta_{top}, \scrSF) \otimes
\Bbb Q) = 0$. This proves (b).
It also follows from \eqref{eq2} that we obtain the isomorphism of
\cite[Proposition~2.1(v), p. 69]{O:CBA}
\begin{equation}
\label{eq10}
H^0(\Delta_{top},\scrSF) \otimes \Bbb Q \cong H^0(\Delta_{top},\scrW)
\otimes \Bbb Q \text{.}
\end{equation}

By \cite[Lemma~8]{DFM:CBg} there is an exact sequence
\begin{equation}
\label{eq6}
0 \to M_0 \to \scrSF(\Delta_{top}) \to \Pic{X} \to 0 \text{.}
\end{equation}
\noindent
Combining \eqref{eq0} and \eqref{eq6},
we have a commutative diagram with exact rows and columns.
\begin{equation}
\label{eq9}
\begin{CD}
0 @. 0 @. 0 \\
@VVV @VVV @VVV \\
M @>>> \scrSF(\Delta_{top}) @>>> \Pic{X} @>>> 0 \\
@VV=V @VVV @VVV \\
M @>>> \scrW(\Delta_{top}) @>>> \Cl(X) @>>> 0
\end{CD}
\end{equation}
Because the center vertical arrow in \eqref{eq9} tensored with $\Bbb Q$ is
the isomorphism \eqref{eq10}, from \eqref{eq9} it follows that
\begin{equation}
\label{eq11}
\Pic(X) \otimes \Bbb Q \cong \Cl(X) \otimes \Bbb Q \text{.}
\end{equation}
It follows from Theorem~\ref{th1} that $\rho_1 = \#(\Delta(1))-s$.
This proves (a). In case (c), $s=r$ so $\rho_1$ is a topological invariant.

(d)
{F}rom Theorem~\ref{th1} and parts (a) and (b),
\[
\rho_0-\rho_1+\rho_2 =
 (r - s) - (\#(\Delta(1)) -s) + 0 = r  - \#(\Delta(1))
\]
which only depends on $\Delta_{top}$.
\end{pf}

\begin{lemma}
\label{lem2}
For any cone $\sigma \in \Delta$, let $\Delta(\sigma)$ denote the subfan of
$\Delta$ consisting of $\sigma$ and all of its faces and $U_\sigma =
\TNemb(\Delta(\sigma))$.
Then $H^0(\Delta(\sigma)_{top},\scrP) = \Cl(U_\sigma)$.
\end{lemma}
\begin{pf}
For each $\sigma \in \Delta$ we have $H^1(\Delta(\sigma)_{top}, \scrSF)
= 0$ \cite[Lemma 2.a, p. 139]{DFM:CBg} so from \eqref{eq2}
$$H^0(\Delta(\sigma)_{top},\scrP) = \scrW(\Delta(\sigma)_{top}) /
\scrSF(\Delta(\sigma)_{top}) \text{.}$$
Now $\scrW(\Delta(\sigma)_{top}) = \Bbb Z^{\#(\Delta(\sigma)(1))}$ and
support functions are linear on a cone $\sigma$, so
$$\scrW(\Delta(\sigma)_{top}) /  \scrSF(\Delta(\sigma)_{top})
\cong \Bbb Z^{\#(\Delta(\sigma)(1))}/ \im(M)
= \Cl(U_\sigma) \text{.}$$
\end{pf}

\begin{lemma}
\label{lem3}
Let $\sigma$ be a cone in $N \otimes \Bbb R$ and $s = \dim{\sigma}$. Then
$\dim_{\Bbb Q}(\Cl(U_\sigma) \otimes \Bbb Q) = \#(\Delta(\sigma)(1))-s$.
Also $\sigma$ is simplicial if and only if  $\Cl(U_\sigma) $ is torsion.
\end{lemma}
\begin{pf}
Follows from Theorem~\ref{th1}.
\end{pf}

The following can be considered a theorem for 3-dimensional fans.

\begin{theorem}
\label{th3}
Let $\Delta$ be a fan on $N \otimes \Bbb R$.
Let $\sigma_0, \dots, \sigma_w$ be the maximal cones in $\Delta$.
Assume
$\sigma_i \cap
\sigma_j$ is simplicial for each $i \not = j$. These assumptions are
satisfied for example if $\dim{\Delta_{top}} \le 3$.
Then
\begin{itemize}
\item[(a)]
\[
 \rho_1 +s + \sum_{i=0}^w(\#(\Delta(\sigma_i)(1))-s_i) =
 \rho_2 +  \#( \Delta(1)) \text{,}
\]
where  we set $s_i = \dim{\sigma_i}$ for each $i = 0, \dots, w$ and $s =
\dim_{\Bbb R} \Bbb R |\Delta_{fan}|$.
\item[(b)]
$\rho_0-\rho_1+\rho_2$ is a topological invariant of $\Delta_{fan}$.
\end{itemize}
\end{theorem}

\begin{pf} (a)
The set $\{ \Delta(\sigma_i)_{top}\}_{i=0}^w$ is an open cover of
$\Delta_{top}$ and the sequence
\begin{equation}
\label{eq3}
0 \to H^0(\Delta_{top},\scrP) \to \bigoplus_{i=0}^w
H^0(\Delta(\sigma_i)_{top},\scr P) \to
\bigoplus_{i=1}^w\bigoplus_{j=0}^{i-1}H^0(\Delta(\sigma_i \cap
\sigma_j)_{top},\scrP)
\end{equation}
is exact since $\scrP$ is a sheaf. Applying Lemma~\ref{lem2}, the sequence
\eqref{eq3} can be written
\begin{equation}
\label{eq4}
0 \to H^0(\Delta_{top},\scrP) \to \bigoplus_{i=0}^w \Cl(U_{\sigma_i}) \to
\bigoplus_{i=1}^w\bigoplus_{j=0}^{i-1} \Cl(U_{\sigma_i \cap \sigma_j}) \text{.}
\end{equation}
By our assumption $\sigma_i \cap \sigma_j$  is
simplicial. By Lemma~\ref{lem3}, $\Cl(U_{\sigma_i \cap \sigma_j})$ is
torsion. By \cite[Theorem 1]{DFM:CBg} $H^1(\Delta_{top}, \scrSF) \cong
H^2(K/X_{\et},\Bbb G_m)$
and the torsion-free part of $H^2(X_{\et},\Bbb G_m)$ is equal
to the torsion-free part of $H^2(K/X_{\et},\Bbb G_m)$. We compute the rank of
the
torsion-free part of $H^2(K/X_{\et},\Bbb G_m)$ from \eqref{eq2} tensored with
$\Bbb Q$:
\begin{equation}
\label{eq5}
0 \to \scrSF(\Delta_{top}) \otimes \Bbb Q \to \scrW(\Delta_{top})  \otimes
\Bbb Q \to
\scrP(\Delta_{top}) \otimes \Bbb Q \to H^2(K/X_{\et}, \Bbb G_m) \otimes
\Bbb Q \to 0  \text{.}
\end{equation}

Tensoring \eqref{eq6} with $\Bbb Q$ and counting dimensions we find
$\dim_{\Bbb Q}(\scrSF(\Delta_{top}) \otimes \Bbb Q) =$ $ \dim_{\Bbb
Q}(\Pic{X}$ $\otimes \Bbb Q) + s$ $ = \rho_1 +s$. By definition
$\scrW(\Delta_{top}) \otimes \Bbb Q
= \Bbb Q^{\#(\Delta(1))}$.
{}From \eqref{eq4} $\scrP(\Delta_{top}) \otimes \Bbb Q
= \bigoplus_{i=0}^w \Cl(U_{\sigma_i}) \otimes \Bbb Q$.
{F}rom
Lemma~\ref{lem3}, $\dim_{\Bbb Q}(\Cl(U_{\sigma_i}) \otimes \Bbb Q) =
\#(\Delta(\sigma_i)(1))-s_i$.
Counting dimensions in
\eqref{eq5}, we have the equation
\begin{equation}
\label{eq7}
\rho_2 = \rho_1 +s + \sum_{i=0}^w(\#(\Delta(\sigma_i)(1))-s_i) - \#(\Delta(1))
\text{ .}
\end{equation}

(b)
{F}rom (a) and Theorem~\ref{th1} we have
\begin{align*}
\rho_0-\rho_1+\rho_2
= & (r-s) +s +
 \sum_{i=0}^w(\#(\Delta(\sigma_i)(1))-s_i) - \# (\Delta(1)) \\
\mbox{} = & r+  \sum_{i=0}^w(\#(\Delta(\sigma_i)(1))-s_i) - \# (\Delta(1))
\end{align*}
which depends only on $\Delta_{top}$.
\end{pf}

As the next example shows, $\rho_1$ and $\rho_2$ are not topological
invariants of $\Delta_{fan}$ when $r \ge 3$ and $\Delta$ is not simplicial.

\begin{example}
\label{ex1}
Let $\Delta$ be a fan on $\Bbb R^3$ and suppose $\Delta$ consists of
three cones of dimension 3 and 6 cones of dimension 1 such that for each
$\sigma_i \in \Delta(3)$, $\#(\Delta(\sigma_i)(1)) = 4$. Assume that the
intersection of the fan $\Delta$ with the unit sphere traces a graph that
looks like that shown in Figure~\ref{fig1}.

\setlength{\unitlength}{.005in}
\begin{figure}
\center{
\begin{picture}(400,400)
\thinlines
\put(50,50){\line(1,0){300}}
\put(50,50){\line(5,6){250}}
\put(300,350){\line(1,-6){50}}
\put(50,50){\line(2,1){100}}
\put(150,100){\line(1,0){150}}
\put(150,100){\line(1,1){100}}
\put(300,100){\line(1,-1){50}}
\put(250,200){\line(1,-2){50}}
\put(250,200){\line(1,3){50}}

\put(50,50){\circle*{5}}
\put(150,100){\circle*{5}}
\put(250,200){\circle*{5}}
\put(300,100){\circle*{5}}
\put(350,50){\circle*{5}}
\put(300,350){\circle*{5}}

\put(45,50){\makebox(0,0)[r]{\small{5}}}
\put(355,50){\makebox(0,0)[l]{\small{4}}}
\put(145,100){\makebox(0,0)[br]{\small{2}}}
\put(305,100){\makebox(0,0)[bl]{\small{1}}}
\put(245,200){\makebox(0,0)[br]{\small{0}}}
\put(300,355){\makebox(0,0)[b]{\small{3}}}
\put(180,180){\makebox(0,0)[l]{\small{$\sigma_0$}}}
\put(240,75){\makebox(0,0)[r]{\small{$\sigma_2$}}}
\put(310,180){\makebox(0,0)[r]{\small{$\sigma_1$}}}

\end{picture} }
\vspace{-0.35in}
\caption{}
\label{fig1}
\end{figure}

For any such fan $\Delta$, $\Delta_{top}$ is unique up to homeomorphism. We
consider 2 such fans $\Delta$ and $\Delta'$ such that $\rho_1(\Delta) \not
= \rho_1(\Delta')$ and $\rho_2(\Delta) \not
= \rho_2(\Delta')$.

For $\Delta$, take $\Delta(1)$ to be $\{ \Bbb R_{\ge} \eta_i | i=0..5\}$
where $\{\eta_0, \dots, \eta_5\}$ $ =$
\begin{equation} \label{eq23.5}
 \left\{
\begin{pmatrix} 1 \\ 0 \\ -2 \end{pmatrix},
\begin{pmatrix} -1 \\ 2 \\ -2 \end{pmatrix},
\begin{pmatrix} -1 \\ -2 \\ -2 \end{pmatrix},
\begin{pmatrix} 1 \\ 0 \\ 2 \end{pmatrix},
\begin{pmatrix} -1 \\ 2 \\ 2 \end{pmatrix},
\begin{pmatrix} -1 \\ -2 \\ 2 \end{pmatrix} \right\}
\text{.}
\end{equation}
Using the methods of \cite[Section 4]{F:Elt} we find that
$\rho_1(\Delta) =1$ and  $\rho_2(\Delta) =1$.

For $\Delta'$, take $\Delta'(1)$ to be $\{ \Bbb R_{\ge} \eta_i | i=0..5\}$
where $\{ \eta_0, \dots, \eta_5\}$ $ =$
\[ \left\{
\begin{pmatrix} 0 \\ 1 \\ 1 \end{pmatrix},
\begin{pmatrix} 0 \\ 0 \\ 1 \end{pmatrix},
\begin{pmatrix} 1 \\ 0 \\ 1 \end{pmatrix},
\begin{pmatrix} -1 \\ 3 \\ 1\end{pmatrix},
\begin{pmatrix} -2 \\ -1 \\ 1\end{pmatrix},
\begin{pmatrix} 3 \\ -1 \\ 1 \end{pmatrix} \right\} \text{.} \]
Using the methods of \cite[Section 4]{F:Elt} we find that
$\rho_1(\Delta') =0$ and  $\rho_2(\Delta') =0$.
\end{example}

\section{A Stratification of the fibers of $\frak T$}
\label{sec4}
Let $\Delta$ be a fan on $N \otimes \Bbb R= \Bbb R^r$ with $\Delta(1) =
\{r_0, \dotsc, r_n \}$. The intersection of $\Delta(1)$ with the unit
sphere $S$ in $\Bbb R^r$ is a finite set of points, say $\{ p_0, \dotsc,
p_n\}$.
About each $p_i$ we can find an open ball $B_i$ on $S$ such that if
$p_i$ is parametrized by $B_i$, then each choice of $\vec{p} = (p_0, p_1,
\dotsc,
p_n)$ in $B_0 \times B_1 \times \dotsm \times B_n$ defines a fan $\Phi =
\Phi(\vec{p})$ such that $\Phi_{top} \cong \Delta_{top}$. The manifold
$\dis{B = \prod_{i=0}^n B_i}$
parametrizes a  subset of fans in the fiber
$\frak T^{-1}(\Delta_{top})$.  Call $B$ an {\em open neighborhood of
$\Delta$}.
If $\vec{p} \in B$, then the fan $\Phi = \Phi(\vec{p})$ is not necessarily
rational.
Sometimes it will be necessary to refer to
points in $B$ that give rise to rational fans. In this case let
\begin{multline}
\label{eq35}
B_{rat} = \{ (p_0, \dotsc, p_n) | \text{ for each } i \text{, }  \\
p_i \text{ is the
intersection of a rational 1-dimensional cone } r_i \text{ with } B_i \}
\text{.}
\end{multline}
For the present section only
we define the set of support functions on a fan to be a real vector space.
If $\sigma$ is a cone, define $\scrSF(\sigma) $ to be $\Hom_{\Bbb R}(\Bbb R
\sigma, \Bbb R)$. Define $\scrSF(\Delta)$ to be the kernel of $\delta^0$ in
the $\Cech$ complex
\begin{equation}
\label{eq34}
0 \to \underset{i}{\oplus}
\scrSF(\sigma_i) \stackrel{\delta^0}{\rightarrow} \underset{i<j}{\oplus}
\scrSF(\sigma_{ij}) \stackrel{\delta^1}{\rightarrow} \underset{i<j<k}{\oplus}
\scrSF(\sigma_{ijk}) \to \dots
\end{equation}
where $\{ \sigma_0, \dotsc, \sigma_w \}$ is the set of maximal cones of
$\Delta$. Define $\kappa_0(\Delta) = \dim_{\Bbb R} \scrSF(\Delta)$. If
$\Delta$ is a rational fan, then this definition of $\kappa_0$ agrees with
the definition given in Remark~\ref{re7} of Section~\ref{sec3}.
In this section we consider the stratification of the manifold $B$ by the
invariant $\kappa_0$.

\begin{example}
\label{ex15} If $\Delta$ is simplicial, then $\kappa_0 = \#(\Delta(1)) =
n+1$ so $B$ has only  1 stratum. As was suggested in Example~\ref{ex1}, we
expect the stratification to be more interesting when $\Delta$ is
nonsimplicial.
\end{example}

\begin{example}
\label{ex16} Let $\Delta$ be the fan on $\Bbb R^3$ given in
equation \eqref{eq23.5} of Example~\ref{ex1}.
Let $\dis{B = \prod_{i=0}^5 B_i}$ be an open neighborhood of $\Delta$.
One can check that any support function $h \in \scrSF(\Delta)$ is
completely determined by its values on $r_0, r_1, r_2, r_3$ so
$\kappa_0(\Delta) \le 4$. {F}rom Example~\ref{ex1} we know that
$\kappa_0(\Delta) = 4$.
It is possible to vary any one of the $r_i$ to
achieve a fan $\Phi$ in $B$ with $\kappa = 3$.
So $B$ has exactly 2 strata. We will see later that the stratum where
$\kappa_0 = 4$ is a Zariski closed subset of $B$.
\end{example}

\begin{conjecture}
\label{conj1}
Let $\Delta$ be a complete fan on $\Bbb R^3$  such that
for each cone $\sigma \in \Delta(3)$,
$\sigma$ is nonsimplicial. Let $B$ be an open neighborhood of $\Delta$ as
described above.
Then for a general choice of $\vec{p} \in B$, if
$\Phi = \Phi(\vec{p})$, then every $\Phi$-linear support
function is linear. In particular for a general choice of $\vec{p} \in
B_{rat}$,
$\kappa_0(\Phi) = 3$ hence $\rho_1(\Phi) = 0$ and
$\rho_2(\Phi)$  is a topological invariant.
\end{conjecture}

In Conjecture~\ref{conj1} by ``general choice'' of
$\vec{p}$ we mean that there is a dense open subset $G \subseteq B$ and
each fan in the
set  $\{ \Phi(\vec p) | \vec{p} \in G \}$ satisfies
the conjecture. That is, if Conjecture~\ref{conj1} is true, a sufficiently
general fan $\Delta'$ with $\Delta'_{top} \cong \Delta_{top}$
should satisfy $\kappa_0(\Delta') = 3$.

As motivation for Conjecture~\ref{conj1}, consider the case where each
$\sigma \in \Delta(3)$ has exactly 4 1-dimensional faces. Let $\Delta(3) =
\{\sigma_0, \dots, \sigma_w \}$, $\Delta(2) =
\{\tau_0, \dots, \tau_e \}$, $\Delta(1) =
\{r_0, \dots, r_n \}$. The intersections of the cones in $\Delta(2)$ with
the unit sphere $S$ in $\Bbb R^3$ trace out the edges of a graph on $S$.
This graph has $e+1$ edges, $n+1$ vertices and $w+1$ regions. So $w+1 =
(e+1)-(n+1)+2$. Each $\sigma_j$ has exactly 4 $\tau_i$'s and each  $\tau_i$
is in exactly 2 $\sigma_j$'s, so
$2(e+1)=4(w+1)$ or $e+1 = 2(w+1)$. Hence $w+1=n-1$.
{F}rom Theorem~\ref{th1}~(c) $\rho'_1 = \dim_{\Bbb Q}(\Cl(X) \otimes
\Bbb Q) = (n+1)-3 = n-2$ and  $\rho'_1(\Delta(\sigma_i)) = \dim_{\Bbb
Q}(\Cl(U_{\sigma_i}) \otimes \Bbb Q) = 4-3 = 1$. From \eqref{eq5} we have
an exact sequence
\begin{equation}
\label{eq13}
0 \to \scrSF(\Delta_{top}) \otimes \Bbb Q \to \scrW(\Delta_{top}) \otimes
\Bbb Q \to  \bigoplus_{i=0}^w (\Cl(U_{\sigma_i}) \otimes \Bbb Q)
\end{equation}
\noindent
Since linear Weil divisors correspond to linear Cartier divisors
\eqref{eq13} gives rise to
\begin{equation}
\label{eq14}
0 \to \Pic(X) \otimes \Bbb Q \to \Cl(X) \otimes
\Bbb Q \to  \bigoplus_{i=0}^w (\Cl(U_{\sigma_i}) \otimes \Bbb Q)
\end{equation}
\noindent
In \eqref{eq14} the middle term has dimension $n-2$ and the third term
dimension $n-1$. For each $i$ the map $\Cl(X) \to \Cl(U_{\sigma_i})$ is
surjective. So we can view $ \Pic(X) \otimes \Bbb Q$ as the intersection of
$n-1$ hyperplanes through $(0)$ in $\Bbb Q^{n-2}$. In general, this
intersection should be $(0)$.
If the conclusion of Conjecture~\ref{conj1} is satisfied, then from
\eqref{eq7} it follows that $\rho_2 = 3 + (w+1)-(n+1)$ (for a general
choice of $\Delta_{fan}$).

Next we give an algorithm for computing an upper bound for $\kappa_0$ for a
fan of arbitrary dimension.
The algorithm can also be used to obtain an upper bound for $\rho_1$ and
$\rho_2$ for complete rational 3-dimensional fans.
If $\Delta$ is a  fan on $\Bbb R^3$ which contains at least 1 cone of
dimension 3, then $\rho_0 = 0$ and by Theorem~\ref{th3} we have
$\rho_1 = \rho_2 + (\text{topological invariant})$.
In this setting   $\rho_1 = \kappa_0-3$ .

\begin{algorithm}
\label{alg1}
Let $\Delta$ be a fan on $N \otimes \Bbb R$. The following is an
algorithm for computing an upper bound for $\kappa_0$.
\end{algorithm}

The algorithm is based on the fact that the map $\scrSF(\Delta) \to
\Bbb Z^{\# (\Delta(1))}$ is injective.
The algorithm finds a subset $G$ of
$\Delta(1)$ such that any support function $h$ in
$\scrSF(\Delta) \otimes \Bbb Q$ is completely determined by its values on
the 1-dimensional cones in $G$.

If $\sigma$ is a maximal cone in $\Delta$ of dimension $d$, then a support
function $h$ is determined by its values on any $d$ 1-dimensional faces of
$\sigma$ that span a $d$-dimensional subspace of $N\otimes \Bbb R$. Pick
$d$ such elements of $\sigma(1)$ and place them in a set called $G$. Place
all other elements of $\sigma(1)$ in a set called $R$. Initially, $G$ and
$R$ are both empty, and the starting cone $\sigma$ is chosen somewhat
arbitrarily. The algorithm proceeds to branch from the starting cone
$\sigma$ outward until all maximal cones of $\Delta$ have been visited and
$\Delta(1)$ has been partitioned into $\Delta(1)=G\cup R$. The order in
which the maximal cones are traversed is somewhat arbitrary and may affect
both the resulting set $G$ and the resulting cardinality of $G$.

\begin{itemize}
\item[Step 0.]
Set $B= \{ \sigma \in \Delta | \sigma \text{ is a maximal cone in } \Delta \}$.
Set
$G= \emptyset$  and $R= \emptyset$.  Go to Step 3.
\item[Step 1.]
If there is a maximal cone $\sigma \in B$ such that
$\sigma(1) \cap (G \cup R)$
contains a spanning set for $\Bbb R \sigma$,
then add the remaining cones in $\sigma(1)-G-R$ to $R$. Remove $\sigma$
from $B$. repeat Step~1 until the condition is false.
\item[Step 2.]
If there is a maximal cone $\sigma \in B$ such that
$\sigma(1) \cap (G \cup R) \not = \emptyset$, then pick $\sigma \in B$ such
that
\begin{itemize}
\item[(i)] $e = \dim_{\Bbb R} \langle \sigma(1) \cap (G \cup R) \rangle$ is
maximal and
\item[(ii)] $d = \dim{\sigma}$ is maximal among all $\sigma \in B$
satisfying (i).
\end{itemize}
For any $\sigma$ satisfying (i) and (ii), pick $\tau_1, \dots, \tau_e$ in
$\sigma(1) \cap (G \cup R)$ such that $\tau_1+ \dots+ \tau_e$ has dimension
$e$. Choose
$\tau_{e+1}, \dots, \tau_d$ in $\sigma(1)$ such that  $\tau_1+ \dots+
\tau_d$ has dimension $d$. Add $\tau_{e+1}, \dots, \tau_d$ to $G$ and add
the remaining elements $\sigma(1)-G-R-\{ \tau_{e+1}, \dots, \tau_d \}$ to $R$.
Delete $\sigma$ from $B$. Go to Step 1.
\item[Step 3.]
If $B \not = \emptyset$, then pick $\sigma \in B$ such that $d =
\dim{\sigma}$ is maximal. Pick
$\tau_{1}, \dots, \tau_d$
in $\sigma(1)$ such that  $\tau_1+ \dots+ \tau_d$
has dimension $d$. Add $\tau_{1}, \dots, \tau_d$ to $G$ and add the
remaining cones in  $\sigma(1)-\{ \tau_{1}, \dots, \tau_d \}$ to $R$.
Delete $\sigma$ from $B$. Go to Step 1.
\item[Step 4.] This point is reached only if $B = \emptyset$. Now
$\Delta(1)$ is partitioned into 2 sets: $\Delta(1) = G \cup R$. Any support
function $h$ in $\scrSF(\Delta) \otimes \Bbb Q$ is determined completely by
its values on $G$. So $\scrSF(\Delta) \to \Bbb Z^{\#(G)}$ is injective.
Therefore $\#(G)$ is an upper bound for $\kappa_0$.
\end{itemize}

\begin{example}
\label{ex5}
Let $\Delta$ be a fan on $\Bbb R^3$ that
consists of three 3-dimensional cones and assume that the
intersection of the fan $\Delta$ with the unit sphere $S$ traces a graph
as shown in Figure~\ref{fig2}(a).
In this example, we step through Algorithm~\ref{alg1} to see that
$\kappa_0(\Delta) \le 4$.
It is shown later in Example~\ref{ex2} that for this fan, $\kappa_0=4$.
Initially, $B=\{\sigma_0,\sigma_1,\sigma_2\}$ and $G=R=\emptyset$.
The algorithm proceeds to Step~3. Place $r_1$, $r_4$,
$r_2$ from $\sigma_0(1)$ in $G$ and $r_0$ in $R$. Delete $\sigma_0$ from
$B$. The condition in Step~1 is still false, so the algorithm goes to
Step~2. For $\sigma_1$, $r_0$ and $r_2$ are both in $G\cup R$ and  $r_0+r_2$
has dimension $e=2$. Place $r_5$ in $G$ and $r_3$ in $R$. Delete $\sigma_1$
from $B$. Go to Step~1. This time the set $G\cup R$ contains $\{ r_0, r_1,
r_3\}$ which is a spanning set for $\Bbb R \sigma_2$. Therefore,
remove $\sigma_2$  from $B$ and
place $r_6$  in $R$. Any support function $h$ is completely determined by
its values on $r_1$, $r_2$, $r_4$ and $r_5$, so $\kappa_0 \le 4$.
\end{example}
\setlength{\unitlength}{.005in}
\begin{figure}
\center{
\hfill (a)
\begin{picture}(400,400)
\thinlines
\put(200,200){\line(0,1){150}}
\put(75,125){\line(5,3){125}}
\put(200,200){\line(5,-3){125}}
\put(75,125){\line(0,1){150}}
\put(75,125){\line(5,-3){125}}
\put(75,275){\line(5,3){125}}
\put(200,50){\line(5,3){125}}
\put(200,350){\line(5,-3){125}}
\put(325,125){\line(0,1){150}}

\put(200,200){\circle*{5}}
\put(200,350){\circle*{5}}
\put(325,125){\circle*{5}}
\put(75,125){\circle*{5}}
\put(325,275){\circle*{5}}
\put(200,50){\circle*{5}}
\put(75,275){\circle*{5}}

\put(195,205){\makebox(0,0)[br]{\small{$r_0$}}}
\put(200,355){\makebox(0,0)[b]{\small{$r_1$}}}
\put(330,125){\makebox(0,0)[l]{\small{$r_2$}}}
\put(70,125){\makebox(0,0)[r]{\small{$r_3$}}}
\put(330,275){\makebox(0,0)[l]{\small{$r_4$}}}
\put(200,45){\makebox(0,0)[t]{\small{$r_5$}}}
\put(70,275){\makebox(0,0)[r]{\small{$r_6$}}}

\put(250,250){\makebox(0,0)[c]{\small{$\sigma_0$}}}
\put(200,125){\makebox(0,0)[c]{\small{$\sigma_1$}}}
\put(150,250){\makebox(0,0)[c]{\small{$\sigma_2$}}}
\end{picture}
\hfill (b)
\setlength{\unitlength}{.005in}
\begin{picture}(400,400)
\thinlines
\put(50,50){\line(1,0){300}}
\put(50,50){\line(1,1){100}}
\put(50,50){\line(0,1){300}}
\put(50,350){\line(1,0){300}}
\put(50,350){\line(1,-1){100}}
\put(150,150){\line(1,0){100}}
\put(150,150){\line(0,1){100}}
\put(150,250){\line(1,0){100}}
\put(250,150){\line(0,1){100}}
\put(250,150){\line(1,-1){100}}
\put(350,50){\line(0,1){300}}
\put(250,250){\line(1,1){100}}

\put(50,50){\circle*{5}}
\put(50,350){\circle*{5}}
\put(150,150){\circle*{5}}
\put(150,250){\circle*{5}}
\put(250,150){\circle*{5}}
\put(250,250){\circle*{5}}
\put(350,50){\circle*{5}}
\put(350,350){\circle*{5}}

\put(250,255){\makebox(0,0)[b]{\small{0}}}
\put(250,145){\makebox(0,0)[t]{\small{1}}}
\put(150,145){\makebox(0,0)[t]{\small{2}}}
\put(150,255){\makebox(0,0)[b]{\small{3}}}
\put(355,355){\makebox(0,0)[bl]{\small{4}}}
\put(355,45){\makebox(0,0)[tl]{\small{5}}}
\put(45,45){\makebox(0,0)[tr]{\small{6}}}
\put(45,355){\makebox(0,0)[br]{\small{7}}}
\put(200,200){\makebox(0,0)[c]{\small{$\sigma_0$}}}
\put(300,200){\makebox(0,0)[c]{\small{$\sigma_1$}}}
\put(200,100){\makebox(0,0)[c]{\small{$\sigma_2$}}}
\put(100,200){\makebox(0,0)[c]{\small{$\sigma_3$}}}
\put(200,300){\makebox(0,0)[c]{\small{$\sigma_4$}}}
\put(10,200){\makebox(0,0)[c]{\small{$\sigma_5$}}}

\end{picture} \hfill 
}
\caption{}
\label{fig2}
\end{figure}

\begin{example}
\label{ex2}
Let $\Delta$ be a complete  fan on $\Bbb R^3$ and assume that the
intersection of the fan $\Delta$ with the unit sphere $S$ traces a graph
that corresponds to the edges of a cube as shown in Figure~\ref{fig2}(b).
Applying Algorithm~\ref{alg1} to $\Delta$, we see that $\rho_1(\Delta) \le
1$ and  $\rho_2(\Delta) \le 1$.

Take $\Delta(1)$ to be $\{ \Bbb R_{\ge} \eta_i | i=0..7\}$
where $\{ \eta_0, \dots, \eta_7\}$ $ =$
\[ \left\{
\begin{pmatrix} 1 \\ 1 \\ 1 \end{pmatrix},
\begin{pmatrix} 1 \\ -1 \\ 1\end{pmatrix},
\begin{pmatrix} -1 \\ -1 \\ 1 \end{pmatrix},
\begin{pmatrix} -1 \\ 1 \\ 1 \end{pmatrix},
\begin{pmatrix} 1 \\ 1 \\ -1 \end{pmatrix},
\begin{pmatrix} 1 \\ -1 \\ -1\end{pmatrix},
\begin{pmatrix} -1 \\ -1 \\ -1\end{pmatrix},
\begin{pmatrix} -1 \\ 1 \\ -1 \end{pmatrix}
\right\} \text{.} \]
Using the methods of \cite[Section 4]{F:Elt} we find that the upper bounds
predicted by Algorithm~\ref{alg1} are reached:
$\rho_1(\Delta) =1$ and  $\rho_2(\Delta) =2$.

Now change the fan so that $\Delta'(1)$ is no longer symmetrical about the
origin. For example,
take $\Delta'(1)$ to be $\{ \Bbb R_{\ge} \eta_i | i=0..7\}$
where $\{ \eta_0, \dots, \eta_7\}$ $ =$
\[ \left\{
\begin{pmatrix} 2 \\ 1 \\ 1 \end{pmatrix},
\begin{pmatrix} 1 \\ -1 \\ 1\end{pmatrix},
\begin{pmatrix} -1 \\ -1 \\ 1\end{pmatrix},
\begin{pmatrix} -1 \\ 1 \\ 1 \end{pmatrix},
\begin{pmatrix} 1 \\ 1 \\ -1 \end{pmatrix},
\begin{pmatrix} 1 \\ -1 \\ -1\end{pmatrix},
\begin{pmatrix} -1 \\ -1 \\ -1\end{pmatrix},
\begin{pmatrix} -1 \\ 1 \\ -1 \end{pmatrix}
\right\} \text{.} \]
Using the methods of \cite[Section 4]{F:Elt} we find that the lower bounds
predicted by Conjecture~\ref{conj1} are attained:
$\rho_1(\Delta') =0$ and  $\rho_2(\Delta') =1$.
Now $\Pic{X'}$ is torsion-free for the complete toric variety $X' =
\TNemb{\Delta'}$. Since $\rho_1(\Delta')=0$ we see that $\Pic{X'}=0$. This
proves that $X'$ is nonprojective (see Remark~\ref{re4}).
If $B$ is an open neighborhood of  $\Delta$, then the strata of $B$ are
$\kappa_0 = 4$ and $\kappa_0 = 3$. We will show later that $\kappa_0 = 4$
corresponds to a Zariski closed subset of $B$.
\end{example}
\begin{remark}
\label{re4}
In general any toric variety satisfying
Conjecture~\ref{conj1} is nonprojective. This is because a projective
normal
variety $X$ will always have a nonprincipal Cartier divisor corresponding
to a hyperplane section.
\samepage{
This follows from commutative
diagram \eqref{eq15}. See \cite[Ex. 6.2, p. 146]{H:AG}.
}
\begin{equation}
\label{eq15}
\begin{CD}
\Pic{ \Bbb P^N} = \Cl(\Bbb P^N) @>>> \Pic{X}  @>>> \Cl(X) \\
 @V{\deg}V{\cong}V @.  @VV{\deg}V \\
  \Bbb Z @>{\cdot(\deg{X})}>> \Bbb Z @= \Bbb Z \\
\end{CD}
\end{equation}
\end{remark}
\begin{example}
\label{ex3}
We give an example to illustrate how \eqref{eq13} can be used to compute
$\kappa_0$. Say $\Delta$ consists of three 3-dimensional cones as shown in
Figure~\ref{fig2}(a).
Then $\scrW(\Delta_{top}) = \Bbb Zr_0 \oplus \dots \oplus \Bbb Zr_6$ and
$\scrW(\Delta(\sigma_0)) = \Bbb Zr_0 \oplus \Bbb Zr_1 \oplus \Bbb Zr_2
\oplus \Bbb Zr_4$.
Let $\eta_i$ be a primitive lattice point on $r_i$, so that $r_i = \Bbb
R_{\ge} \eta_i$ for $i=0..6$.
The kernel of the surjection $\phi_0 : \scrW(\Delta)
\to \Cl(U_{\sigma_0})$ is spanned by the vectors
$ \begin{pmatrix} 0 & 0 & 0 & 1 & 0 & 0 & 0 \end{pmatrix} ^\top$,
$ \begin{pmatrix} 0 & 0 & 0 & 0 & 0 & 1 & 0  \end{pmatrix} ^\top$,
$ \begin{pmatrix} 0 & 0 & 0 & 0 & 0 & 0 & 1  \end{pmatrix} ^\top$,
and the columns of
\(
\begin{pmatrix} \eta_0 & \eta_1 & \eta_2 & 0 & \eta_4 & 0 & 0 \end{pmatrix}
^\top \).
So $\ker{\phi_0}$ is a subspace of codimension 1. Consider the matrix
equation
\begin{equation}
\label{eq16}
\begin{pmatrix}
\eta_0 & \eta_1 & \eta_2 & \eta_4
\end{pmatrix}
 \overrightarrow{v_0} = 0 \text{.}
\end{equation}
Set $A = \begin{pmatrix}
\eta_0 & \eta_1 & \eta_2
\end{pmatrix}$. Then \eqref{eq16}
becomes
\begin{equation}
\label{eq16.6}
\begin{pmatrix}
I & A^{-1} \eta_4
\end{pmatrix}
\overrightarrow{v_0} = 0 \text{.}
\end{equation}
Set $ A^{-1} \eta_4 = \begin{pmatrix}
 a_0 & b_0 & c_0 \end{pmatrix}^ \top $.
So $\overrightarrow{v_0} = \begin{pmatrix} -a_0 z & -b_0 z & -c_0 z & z
\end{pmatrix} ^ \top $.
Since any 3 columns of the matrix in \eqref{eq16} are linearly independent,
$\overrightarrow{v_0}$ has 4
nonzero entries, or $\overrightarrow{v_0}=0$.
Normalize $\overrightarrow{v_0}$ by taking $z = -1$.
Then $\ker{\phi_0}$ is the set of solutions to
\begin{equation}
\label{eq17}
\begin{pmatrix}
a_{0} & b_{0} & c_{0} & 0 & -1 & 0 & 0
\end{pmatrix}
\vec{x} = 0 \text{.}
\end{equation}
Hence if $\phi : \scrW(\Delta) \to \Cl(U_{\sigma_0}) \oplus
\Cl(U_{\sigma_1}) \oplus \Cl(U_{\sigma_2})$, then $\ker{\phi}$ is the set of
solutions to
\begin{equation}
\label{eq18}
\begin{pmatrix}
a_{0} & b_{0} & c_{0} & 0 & -1 & 0 & 0 \\
a_{1} & 0  & b_{1} & c_{1} & 0 & -1 &  0 \\
a_{2} & b_{2} & 0 &  c_{2} & 0 & 0 & -1
\end{pmatrix} \vec{x} = 0 \text{.}
\end{equation}
This coefficient matrix has rank 3 so $\ker{\phi}$ has rank 4. Therefore
$\kappa_0 =4$ and for any open neighborhood $B$ of $\Delta$, $B$ has only 1
stratum.  In this case, we see that $\kappa_0$ and hence $\rho_1$ and
$\rho_2$ are topological invariants of the set of all (rational) fans that
look like the one shown in Figure~\ref{fig2}(a). We could assume
$\Delta$ has more than three (say $w+1$) 3-dimensional cones each with four
1-dimensional faces meeting around the common
1-dimensional face $r_0$. By a similar argument we see that $\kappa_0 = w+2$.
\end{example}

\begin{example}
\label{ex4}
Let $\Delta$ be a fan on $\Bbb R^3$ such that $\Delta_{top}$ is
homeomorphic to the fan in Example~\ref{ex1}.
Following the procedure of Example~\ref{ex3}, set up equations analogous to
\eqref{eq16} \eqref{eq17} and \eqref{eq18}. Then $\ker{\phi}$ is the set of
solutions to
\begin{equation}
\label{eq19}
\begin{pmatrix}
a_{0} & 0 & b_{0} & c_{0} & 0 & -1 \\
a_{1} & b_{1} & 0 & c_{1} & -1 &  0 \\
0 & a_{2} &  b_{2} & 0 & c_{2} & -1
\end{pmatrix}
\vec{x} = 0 \text{.}
\end{equation}
The coefficient matrix in \eqref{eq19} clearly has rank 2 or more. This
agrees with the upper bound 4 predicted for $\kappa_0$ by
Algorithm~\ref{alg1}. The third, fourth and sixth columns of \eqref{eq19}
are independent if and only if
\begin{equation}
\label{eq20}
(b_2-b_0) c_1 \not = 0 \text{ .}
\end{equation}
This shows that on the complement of a Zariski open subset of $B$,
$\kappa_0 = 3$. We check that \eqref{eq20} is satisfied on a nonempty
subset of $B$.
Note that \eqref{eq20} is satisfied if
\begin{equation}
\label{eq20.5}
\text{the second row of }
\begin{pmatrix}
\eta_0 & \eta_2 & \eta_3
\end{pmatrix}^{-1} \eta_5
\not =
\text{the second row of }
\begin{pmatrix}
\eta_1 & \eta_2 & \eta_4
\end{pmatrix}^{-1} \eta_5
\end{equation}
which will be true for a sufficiently general choice of the fan. To see
this, consider letting $p_0$ vary in $B_0$. Then in \eqref{eq20.5} the
matrix
$\begin{pmatrix}
\eta_0 & \eta_2 & \eta_3
\end{pmatrix}^{-1} $ varies but
the matrix
$\begin{pmatrix}
\eta_1 & \eta_2 & \eta_4
\end{pmatrix}^{-1}$ remains constant.

\end{example}

Consider \eqref{eq13} once again. Let $\Delta$ be a complete  fan on $\Bbb
R^3$. Let $B$ be an open neighborhood of $\Delta$.
Proceed as in Examples \ref{ex3} and \ref{ex4}.
Set up the matrix equation $\Phi \vec{x} = 0$ for  $\ker{\phi}$.
Since $M
\to \SF(\Delta'_{top})$ is injective, $\ker{\phi}$ has rank at least 3.
Consider an arbitrary
$(n-2)- \text{by}-(n-2)$ submatrix
$\Phi_0$ of $\Phi$. Then $\Phi_0$ has rank $n-2$
exactly when $\det(\Phi_0) \not = 0$. As in \eqref{eq19} and \eqref{eq20},
we can show that $\det(\Phi_0)=0$ is an equation in no more than $3(n-2)$
variables which are parametrized by points in $B$.
The equation $\det(\Phi_0) =0$
defines a Zariski closed subset of $B$. On the complement of this closed
set $\det(\Phi_0) \not = 0$, $\rank{\Phi} = n-2$ and $\ker{\phi}$ has rank
3. If there is at least one choice of $\Delta_{fan}$ for which
$\det(\Phi_0) \not =0$, then the open set making up the complement of the
determinant equations  will be nonempty, hence the conclusion of
Conjecture~\ref{conj1} will
be satisfied. This shows for example that the general fan which is
topologically homeomorphic to that of Figure~\ref{fig2}(b) satisfies the
conclusion to Conjecture~\ref{conj1}, because in Example~\ref{ex2} an
example is given which shows the
determinants are nonzero on a nonempty Zariski open in $B$.


\end{document}